\documentclass{amsart}
\usepackage{amsfonts}
\input{epsf}
\usepackage{amsmath, amsthm, amssymb}
\usepackage{latexsym}

\theoremstyle{definition}

\theoremstyle{remark}

\numberwithin{equation}{section}

\begin{document}
\title{Hybrid dynamics for currency modeling}

\author{Ted Theodosopoulos}
\author{Alex Trifunovic}
\address{IKOS Research \\ 9 Castle Square, Brighton, East Sussex, BN1 1DZ \\ United Kingdom}
\email{ptaetheo@earthlink.net}
\email{alex@ikosresearch.com}


\date{May 15, 2006.}



\begin{abstract}
We present a simple hybrid dynamical model as a tool to investigate behavioral strategies based on trend following.  
The multiplicative symbolic dynamics are generated using a lognormal diffusion model for the at-the-money 
implied volatility term structure.  Thus, are model exploits information from derivative markets to 
obtain qualititative properties of the return distribution for the underlier.  We apply our model to the JPY-USD exchange rate 
and the corresponding 1mo., 3mo., 6mo. and 1yr. implied volatilities.  Our results indicate that the modulation of autoregressive 
trend following using derivative-based signals significantly improves the fit to the distribution of times between successive 
sign flips in the underlier time series.
\end{abstract}

\maketitle

\section{Introduction}\label{sec:int}

The spot market for major currencies is by far the most liquid financial market in the world.  It is also notorious for its speculative nature, 
with over a trillion dollars changing hands daily.  Most of this transaction volume is driven by endogenous `herd' dynamics, either technical 
or behavioral in nature.  Over the past decade, the progressive abundance of publically available data at increasingly higher frequencies have 
made the FX market a testing ground for innovative models that probe where general equilibrium economics left off, beyond the Efficient Markets 
Hypothesis.  On the one hand, empirically derived `stylized facts' have become widely quoted as hallmarks of this new frontier.  These include 
statistically challenging estimates of nonlinear, dynamics behavior, like clustered volatility (i.e. significantly positive serially correlation 
of realized volatility), long memory (i.e. subexponential decay of return autocorrelation function) and fat tails (i.e. probablity of extreme 
returns far exceeding that implied by the normal) \cite{daco, tak1}.

On the other hand, modelers have been moving their focus from models based on homogeneous, autonomous, sophisticated optimizing rational agents 
to heterogeneous, strongly coupled agents with explicitly limited rationality.  This represents a fundamental paradigm shift, necessitating new 
analytical tools to answer innovative questions.  One modeling direction involves agents that choose among a small set of standard strategies.  
The most commonly studied trade-off is between trend following and contrarian trading.  In this context, the question is how to optimally time 
the agents' strategic switching \cite{homm, tak3, tak5, bouc}.

Various attempts have been made to modulate pure trend following \cite{kear, tak4}.  Our approach extends earlier work in using multiplicative signals \cite{sol1} to 
adjust endogenously the degree and timing of trend following.  In this paper we present a hybrid dynamical system in which the continuous trend 
following process is multiplied by a volatility-driven signal \cite{sol2}.  We exploit derivative market information to construct Boolean logic rules that 
instantiate the desired switching between trend-following and contrarian strategies \cite{tak2, tak6}.

Our process consists of an ${\rm AR} (n)$ model multiplied by the symbolic dynamics of the term structure of at-the-money (ATM) implied 
volatility, which is assumed to diffuse lognormally.  This simple model replicates the fat tails of the return distribution, the positively 
autocorrelated realized volatility and and the geometric phase plot pattern which characterise the actual data but are often absent in pure 
trend following simulations.  Moreover, we are able to show a particularly accurate fit to the right tail of the `inter-slip' times, i.e. the 
times between successive strategic switches from a follower to a contrarian rule.

\section{The Hybrid Model}\label{sec:hybrid}

Let $X_k$ represent the price return time series, $Y_{j,k}$ the ATM volatility for the $j^{\rm th}$ term at time $k$ and $\epsilon_k$ an iid 
standard normal sequence.  Let $h: (-\infty, \infty) \times \left[ \left. 0, \infty \right) \right. \longrightarrow \{-1,0,1\}$ be the Boolean 
function that encodes the desired symbolic dynamics.  Here we will use $h(x,y) = \chi_{(y,\infty)} (x) - \chi_{(-\infty,-y)} (x)$.  Consider 
the three parameter family of hybrid systems given by:
\begin{equation}
X_k = n^{-1} g_k \sum_{i=1}^n X_{k-i} + \alpha \epsilon_k  \label{eq:hybrid}
\end{equation}
for $k>n$, where $\alpha$ is a control parameter and
\begin{equation}
g_k = h \left( \sum_{j=1}^m h \left( Y_{j,k} - Y_{j,k-1}, 0 \right), 1 \right).  \label{eq:symbolic} 
\end{equation}
We will model the volatility process as a driftless $m$-dimensional correlated random walk, i.e. there exists a positive definite, symmetric matrix $C \in {\mathcal R}^{m^2}$ such that for all $t = (t_1, t_2, \ldots, t_m)'$ and for all $k>0$,
$${\bf E} \left[ \exp \left\{ \sqrt{-1} \sum_{j=1}^m \left( Y_{j,k} -Y_{j,k-1} \right) t_j \right\} \right] = \exp \left\{ -{\frac {1}{2}} t'Ct \right\}$$
and for $k_1 \ne k_2$, 
$${\bf E} \left[ \left({\bf Y}_{\cdot,k_1} -{\bf Y}_{\cdot,k_1-1} \right) \left({\bf Y}_{\cdot,k_2} -{\bf Y}_{\cdot,k_2-1} \right) \right] =0.$$

\begin{figure}
\epsfxsize=5in
\epsfbox{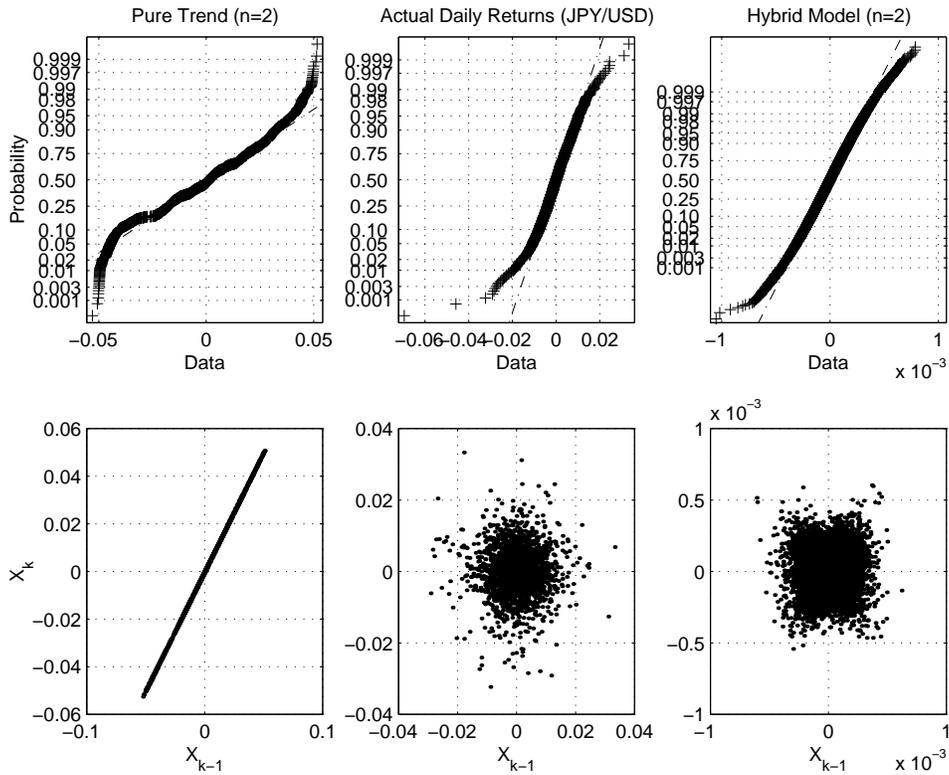}
\caption{The hybrid model matches better the distributional properties of the actual time series than pure trend following.}
\label{fig:fatails}
\end{figure}

We have simulated $10^6$ steps of the hybrid model using the historically determined covariance matrix for the ATM volatility term structure $C$ \cite{jams}.  We have used $m=4$ for all simulations in this paper.  Figure \ref{fig:fatails} shows the comparison of the hybrid model (\ref{eq:hybrid}) on the far right column, against the pure trend following model\footnote{This is equivalent to (\ref{eq:hybrid}) with $g_k \equiv 1$.} in the far left column and the actual daily returns of the JPY/USD exchange rate for five years in the middle column ($2171$ data points).  The first row shows the comparison of the normal probability plots.  This comparison shows that pure trend following generates thin tails, while the hybrid model matches well the fat tails characterizing the actual return distribution.  In both cases $n=2$ was used to generate the graph.

The bottom row shows the comparison of phase plots, i.e. the relationship between $X_k$ and $X_{k-1}$.  Naturally the pure trend following model exhibits very highly serially correlated returns, particularly because we've kept the noise level very low ($\alpha = 10^{-4}$).  Observe that the actual phase plot deviates noticeably from the circular pattern expected from a random walk.  The hybrid model generates a more square-like phase plot pattern, that nevertheless matches the phase plot for the actual data much better than the linear pattern of the pure trend following model.

\begin{figure}
\epsfxsize=4in
\epsfbox{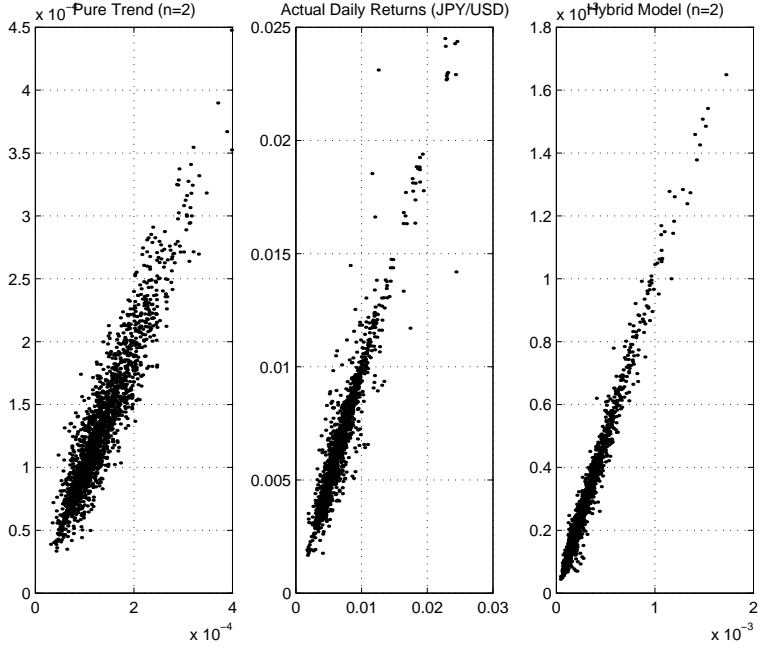}
\caption{Volatility clustering can be seen from the robust positive slope linear patterns in the phase plots of realized volatility.}
\label{fig:clusteredvoly}
\end{figure}

We have also investigated the degree to which the hybrid model generates realized volatility patterns consistent with the clustering phenomenon document in numerous empirical studies.  Figure \ref{fig:clusteredvoly} shows the three way comparison of the biweekly realized volatility phase plot.  In all cases there is a strong positive correlation, which is consistent with the previously reported volatility clustering.  Indeed the resulting pattern from the hybrid model simulation better matches the low dispersion observed in the actual data.

\section{Symbolic Dynamics}\label{sec:symbolic}

We proceed now to investigate in more detail the effect of the symbolic dynamics.  Clearly, the multiplicative modulation of the $AR(n)$ model with $g_k$ serves to flip the sign of the simulated returns at random intervals, determined by the direction of movement of the implied volatility.  
 
Let ${\bf P}^{\alpha,n,m}$ denote the path measure induced by (\ref{eq:hybrid}) with parameters $\alpha$, $n$ and $m$.  Of course
$$\sum_{j=1}^m h \left( Y_{j,k} - Y_{j,k-1}, 0 \right) \in \{-m, -m+2, -m+4, \ldots, m-4, m-2, m\}$$
and therefore
\begin{equation}
q_1 \doteq \lim_{\alpha \rightarrow 0} {\bf P}^{\alpha,1,m} \left(X_{k+1} X_k <0 \right) = {\bf Pr} \left(g_{k+1} <-1 \right)
\label{eq:qonequ}
\end{equation}
which is independent of $k$ since the volatility process has independent increments.  Similarly for $n=2$ we have:
\begin{equation}
q_2 \doteq \lim_{\alpha \rightarrow 0} {\bf P}^{\alpha,2,m} \left( \left\{X_{k+1} (X_k +X_{k-1}) <0 \right\} \cap \left\{ X_k X_{k-1} >0 \right\} \right) = {\bf Pr} \left(g_{k+1} <-1 \right)
\label{eq:qtwoequ}
\end{equation}
Thus, in general,
\begin{equation}
q_n \doteq \lim_{\alpha \rightarrow 0} {\bf P}^{\alpha,n,m} \left( \left\{X_k \sum_{i=1}^n X_{k-i} <0 \right\} \cap \left\{ \bigcap_{i=1}^{n-1} \left\{ X_{k-i} X_{k-i-1} >0 \right\} \right\} \right) = {\bf Pr} \left(g_k <-1 \right)
\label{eq:qnequ}
\end{equation}
Thus $q_1 = q_2 = \ldots = q_n = q$.  The driftless nature of the volatility diffusion implies that the resulting measure on $\{-m, -m+2, -m+4, \ldots, m-4, m-2, m\}$ is symmetric.  Thus, only $ \lceil {\frac {m}{2}} \rceil$ unknowns need to be determined.  Using \ref{eq:qnequ}) we see that
\begin{equation}
q = {\frac {1- {\bf Pr} \left( -1 \leq \sum_{j=1}^m h \left( Y_{j,k} - Y_{j,k-1}, 0 \right) \leq 1 \right)}{2}} \label{eq:qone} 
\end{equation}

\begin{figure}
\epsfxsize=4in
\epsfbox{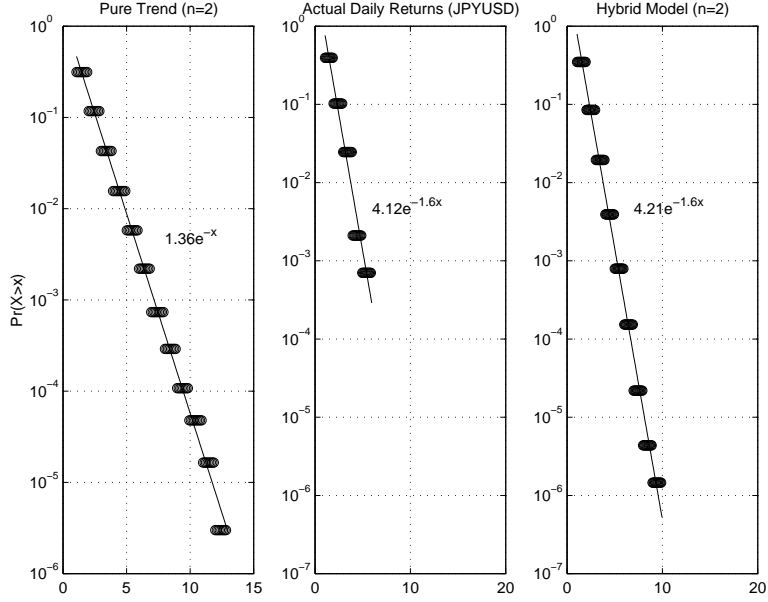}
\caption{The hybrid model matches better the distributional properties of the actual time series than pure trend following.}
\label{fig:interflip}
\end{figure}

Let $T_0=0$ and
$$T_\ell = \min \left\{k> T_{\ell-1} \left| \left\{X_k \sum_{i=1}^n X_{k-i} <0 \right\} \cap \left\{ \bigcap_{i=1}^{n-1} \left\{ X_{k-i} X_{k-i-1} >0 \right\} \right\} \right. \right\}$$
for $\ell>0$ and $\tau_\ell = T_\ell – T_{\ell-1}$.  This sequence of stopping times represents the periods while the process $X_k$ maintains a fixed sign.  The independent increments property discussed above is sufficient in the case $n=1$ to obtain
\begin{equation}
{\bf Pr} \left(\tau_\ell > s \right) = (1 -q)^s. \label{eq:tauone}
\end{equation}
Figure \ref{fig:interflip} shows the complementary cumulative density function of $\tau$.  We observe that the slope of this graph provides a clear indication of the advantage of the hybrid dynamics presented here in modeling the empirical data compared to pure trend following.  Specifically, the simulation of the hybrid model with $n=2$ gives an estimated exponential decay rate of $-1.6$ which is within $3\%$ of the value obtained from the JPY/USD data.  On the other hand, the pure trend following model leads to a much shallower exponential decay rate of $-1$.

Using (\ref{eq:tauone}) we would expect an exponential decay rate equal to $\log(1-q)$ for a hybrid model with $n=1$.  The exponential decay rate for a simulation of the hybrid model with $n=1$ is around $-0.6371$ (data not shown) which would lead to $q \cong 0.4712$.  The resulting measure on $\{-m, -m+2, -m+4, \ldots, m-4, m-2, m\}$ is highly bimodal, as would be expected due to the high positive correlation among the implied volatilities.  

When $n>1$ (\ref{eq:tauone}) no longer holds because of the longer history involved in the rolling average.  Instead of the geometric model which gave rise to (\ref{eq:tauone}) we can use a higher order negative binomial process to accommodated the higher values of $n$.  In particular, for a general $n$, we will approximate the distribution of $\tau_\ell$ with a negative binomial with parameters $q$ and $n$:

\begin{equation}
{\bf Pr} \left(\tau_\ell > s \right) = \left( \begin{array}{c} s-1 \\ n-1 \end{array} \right) (1 -q)^{s-n} q^n. \label{eq:taun}
\end{equation}
The simulation which generated the data on the right-most panel in figure \ref{fig:interflip} used $n=2$.  In this case, using (\ref{eq:taun}) we would expect an exponential decay rate approximately equal to $\log(1-q)-1$, where we have used $\log x \cong x$ for small $x$.  Thus, the observed exponential decay of $-1.6$ corresponds to $q \cong 0.4512$.  This value is slightly lower than that for the case $n=1$ as expected because the averaging of two consecutive terms makes it somewhat more difficult to flip signs (i.e. the probability of long periods with the same sign decays slightly more slowly).

\section{Conclusions and Next Steps}

We have described a simple hybrid dynamical system that models the modulation of trend following strategies for exchange rates using information from the market of currency derivatives.  The resulting family of models was shown to better capture distributional characteristics of the empirical return distribution for the JPY/USD exchange rate.  

The particular form of the symbolic dynamics we used can be interpreted as a qualified majority of the implied volatility signals.  Clearly, different Boolean expressions would lead to different patterns in the distribution of random times $T_\ell$ at which the process changes sign.  It is desirable to classify symbolic dynamics with respect to their effect on the resulting distribution of `inter-flip' times.

So far we have used a driftless model of the volatility term structure.  Adding a true second moment signal (e.g. volatility drift) and an explicit third moment signal (e.g. risk reversal) can serve to make the symbolic modulation of trend following more nuanced.  It would be useful to develop analytical tools, like the exponential decay rate approximation described in section \ref{sec:symbolic}, to detect the effect of different derivative signals in the empirically observable distributional patterns of inter-flip times.

\bibliographystyle{amsalpha}

\begin{thebibliography}{10}
%
\bibitem{homm}
Chiarella, C., He, X.Z., Hommes, C.H.:
A dynamic analysis of moving average rules.
Journal of Economic Dynamics and Control, forthcoming (2006).
%
\bibitem{daco}
Di Matteo, T., Aste, T., Dacorogna, M.M.:
Scaling behaviors in differently developed markets.
Physica A, {\bf 324} (2003), 183-188.
%
\bibitem{jams}
Jamshidian, F., Zhu, Y.:
Scenario simulation: theory and methodology.
Finance and Stochastics, {\bf 1} (1996), 43-67.
%
\bibitem{sol1}
Huang, Z.-F., Solomon, S.:
Stochastic multiplicative processes for financial markets.
Physica A, {\bf 306} (2002), 412-422.
%
\bibitem{kear}
Kearns, M., Kakade, S.:
Trading in Markovian price model.
preprint, \texttt{http://www.cis.upenn.edu/~mkearns/papers/pricemodel.pdf}, in Proceedings of COLT 2005.
%
\bibitem{sol2}
Louzoun, Y., Solomon, S.:
Volatility driven market in a generalized Lotka-Voltera fomralism.
Physica A, {\bf 302} (2001), 220-233.
%
\bibitem{tak1}
Mizuno, T., Kurihara, S., Takayasu, M., Takayasu, H.:
Analysis of high-resolution foreign exchange data of USD-JPY for 13 years.
Physica A, {\bf 324} (2003), 296-302.
%
\bibitem{tak2}
Mizuno T., Nakano, T., Takayasu, M., Takayasu, H.:
Trader's strategy with price feedbacks in financial markets.
Physica A, {\bf 344} (2004), 330-334.
%
\bibitem{tak3}
Mizuno, T., Takayasu, M., Takayasu, H.:
Modeling a foreign exchange rate using moving average of Yen-Dollar market data.
preprint, \texttt{http://www.arxiv.org/abs/physics/0508162}, in Practical Fruits of Econophysics (Springer Verlag Tokyo), 57-61, 2006.
%
\bibitem{tak4}
Ohira, T., Sazuka, N., Marumo, K., Shimizu, T., Takayasu, M., Takayasu, H.:
Predictability of currency market exchange.
Physica A, {\bf 308} (2002), 368-374.
%
\bibitem{tak5}
Ohnishi, T., Mizuno, T., Aihara, K., Takayasu, M., Takayasu, H.:
Statistical properties of the moving average price in dollar-yen exchange rates.
Physica A, {\bf 344} (2004), 207-210.
%
\bibitem{bouc}
Potters, M., Bouchaud, J.-P.:
Trend followers lose more often than they gain, 
preprint, \texttt{http://arxiv.org/pdf/physics/0508104}, 2005.
%
\bibitem{tak6}
Sato, A.-H., Takayasu, H.:
Market price simulator based on analog electric circuit.
preprint, \texttt{http://arxiv.org/pdf/cond-mat/0104318}, 2006.
%
\end{thebibliography}

\end{document}